\documentclass[10.9pt,twoside]{amsart}
\usepackage{amsmath, amsthm, amscd, amsfonts, amssymb, graphicx, color}
\usepackage[bookmarksnumbered, colorlinks, plainpages]{hyperref}

\theoremstyle{definition}

\theoremstyle{remark}

\numberwithin{equation}{section}

\begin{document}
\setcounter{page}{1}
\begin{center}
{\bf  DERIVATIONS AND COHOMOLOGICAL GROUPS OF\\ BANACH ALGEBRAS  }
\end{center}

\title[]{}
\author[]{KAZEM HAGHNEJAD AZAR  }

\address{}

\dedicatory{}



\begin{abstract} Let $B$ be a Banach $A-bimodule$ and let $n\geq 0$. We investigate the relationships between some cohomological groups of $A$, that is, if the topological center of the left module action $\pi_\ell:A\times B\rightarrow B$ of $A^{(2n)}$ on $B^{(2n)}$ is $B^{(2n)}$  and  $H^1(A^{(2n+2)},B^{(2n+2)})=0$, then we have $H^1(A,B^{(2n)})=0$,  and we find the relationships between cohomological groups such as $H^1(A,B^{(n+2)})$  and    $H^1(A,B^{(n)})$,  spacial   $H^1(A,B^*)$ and  $H^1(A,B^{(2n+1)})$. We obtain some results in Connes-amenability of  Banach algebras, and so for  every compact group $G$, we conclude that $H^1_{w^*}(L^\infty(G)^*,L^\infty(G)^{**})=0$.  Let $G$ be an amenable locally compact group. Then there is a Banach $L^1(G)-bimodule$ such as $(L^\infty(G),.)$ such that $Z^1(L^1(G),L^\infty(G))=\{L_{f}:~f\in L^\infty(G)\}.$  We also obtain some conclusions in the Arens regularity of module actions and weak amenability of Banach algebras. We introduce some new concepts as  $left-weak^*-to-weak$ convergence property [$=Lw^*wc-$property] and $right-weak^*-to-weak$ convergence property [$=Rw^*wc-$property] with respect to $A$ and we show that if $A^*$ and $A^{**}$, respectively,  have $Rw^*wc-$property and $Lw^*wc-$property and  $A^{**}$ is weakly amenable, then $A$ is weakly amenable. We also show to relations between a derivation $D:A\rightarrow A^*$ and this new concepts.
\end{abstract} \maketitle

\section{\bf  Preliminaries and
Introduction }

\noindent Let $B$ be a   Banach $A-bimodule$.
   A derivation from $A$ into $B$ is a bounded linear mapping $D:A\rightarrow B$ such that $$D(xy)=xD(y)+D(x)y~~for~~all~~x,~y\in A.$$
The space of continuous derivations from $A$ into $B$ is denoted by $Z^1(A,B)$.\\
Easy example of derivations are the inner derivations, which are given for each $b\in B$ by
$$\delta_b(a)=ab-ba~~for~~all~~a\in A.$$
The space of inner derivations from $A$ into $B$ is denoted by $N^1(A,B)$.
The Banach algebra $A$ is said to be a amenable, when for every Banach $A-bimodule$ $B$, the inner derivations are only derivations existing from $A$ into $B^*$. It is clear that $A$ is amenable if and only if $H^1(A,B^*)=Z^1(A,B^*)/ N^1(A,B^*)=\{0\}$. The concept of amenability for a Banach algebra $A$, introduced by Johnson in 1972, has proved to be of enormous importance in Banach algebra theory, see [13].
A Banach algebra $A$ is said to be a weakly amenable, if every derivation from $A$ into $A^*$ is inner. Similarly, $A$ is weakly amenable if and only if $H^1(A,A^*)=Z^1(A,A^*)/ N^1(A,A^*)=\{0\}$. The concept of weak amenability was first introduced by Bade, Curtis and Dales in [2] for commutative Banach algebras, and was extended to the noncommutative case by Johnson in [14].\\
For Banach $A-bimodule$ $B$, the quotient space $H^1(A,B)$ of all continuous derivations from $A$ into $B$ modulo the subspace of inner derivations is called the first cohomology group of $A$ with coefficients in $B$.\\
In this paper, by using the Arens regularity of module actions, for Banach algebra $A$, we find some relations between cohomology groups $A$ and $A^{(2n)}$ with some applications in the $n-weak$ amenability of Banach algebras that introduced by Dales,  Ghahramani,   Gr{\o}nb{\ae}k in [6]. So for this aim, we extended some propositions from [6, 7, 10] into general situations. We  investigated to relationships between cohomology groups of $A\oplus B$ and $A, B$ where $A$ and $B$ are Banach algebras and on the other hand, we establish some relationships between the Connes-amenability and supper-amenability of Banach algebras with some results in group algebras. We have also some conclusions in the Arens regularity  of module actions. In last section, for Banach $A-module$ $B$, we introduce a new definitions as  $left-weak^*-to-weak$ convergence property [ $=Lw^*wc-$property] and $right-weak^*-to-weak$ convergence property [ $=Rw^*wc-$property] with respect to $A$ and we show that if $A^*$ and $A^{**}$, respectively,  have $Rw^*wc-$property and $Lw^*wc-$property and  $A^{**}$ is weakly amenable, then $A$ is weakly amenable. \\
We introduce some notations and definitions that we used
throughout  this paper.\\ Let $A$ be  a Banach algebra and $A^*$,
$A^{**}$, respectively, are the first and second dual of $A$.  For $a\in A$
 and $a^\prime\in A^*$, we denote by $a^\prime a$
 and $a a^\prime$ respectively, the functionals on $A^*$ defined by $<a^\prime a,b>=<a^\prime,ab>=a^\prime(ab)$ and $<a a^\prime,b>=<a^\prime,ba>=a^\prime(ba)$ for all $b\in A$.
   The Banach algebra $A$ is embedded in its second dual via the identification
 $<a,a^\prime>$ - $<a^\prime,a>$ for every $a\in
A$ and $a^\prime\in
A^*$.  Let $A$ be a Banach algebra. We
say that a  net $(e_{\alpha})_{{\alpha}\in I}$ in $A$ is a left
approximate identity $(=LAI)$ [resp. right
approximate identity $(=RAI)$] if,
 for each $a\in A$,   $e_{\alpha}a\longrightarrow a$ [resp. $ae_{\alpha}\longrightarrow a$].\\
 \noindent Let $X,Y,Z$ be normed spaces and $m:X\times Y\rightarrow Z$ be a bounded bilinear mapping. Arens in [1] offers two natural extensions $m^{***}$ and $m^{t***t}$ of $m$ from $X^{**}\times Y^{**}$ into $Z^{**}$ as following\\
1. $m^*:Z^*\times X\rightarrow Y^*$,~~~~~given by~~~$<m^*(z^\prime,x),y>=<z^\prime, m(x,y)>$ ~where $x\in X$, $y\in Y$, $z^\prime\in Z^*$,\\
 2. $m^{**}:Y^{**}\times Z^{*}\rightarrow X^*$,~~given by $<m^{**}(y^{\prime\prime},z^\prime),x>=<y^{\prime\prime},m^*(z^\prime,x)>$ ~where $x\in X$, $y^{\prime\prime}\in Y^{**}$, $z^\prime\in Z^*$,\\
3. $m^{***}:X^{**}\times Y^{**}\rightarrow Z^{**}$,~ given by~ ~ ~$<m^{***}(x^{\prime\prime},y^{\prime\prime}),z^\prime>$  $=<x^{\prime\prime},m^{**}(y^{\prime\prime},z^\prime)>$\\ ~where ~$x^{\prime\prime}\in X^{**}$, $y^{\prime\prime}\in Y^{**}$, $z^\prime\in Z^*$.\\
The mapping $m^{***}$ is the unique extension of $m$ such that $x^{\prime\prime}\rightarrow m^{***}(x^{\prime\prime},y^{\prime\prime})$ from $X^{**}$ into $Z^{**}$ is $weak^*-to-weak^*$ continuous for every $y^{\prime\prime}\in Y^{**}$, but the mapping $y^{\prime\prime}\rightarrow m^{***}(x^{\prime\prime},y^{\prime\prime})$ is not in general $weak^*-to-weak^*$ continuous from $Y^{**}$ into $Z^{**}$ unless $x^{\prime\prime}\in X$. Hence the first topological center of $m$ may  be defined as following
$$Z_1(m)=\{x^{\prime\prime}\in X^{**}:~~y^{\prime\prime}\rightarrow m^{***}(x^{\prime\prime},y^{\prime\prime})~~is~~weak^*-to-weak^*-continuous\}.$$
Let now $m^t:Y\times X\rightarrow Z$ be the transpose of $m$ defined by $m^t(y,x)=m(x,y)$ for every $x\in X$ and $y\in Y$. Then $m^t$ is a continuous bilinear map from $Y\times X$ to $Z$, and so it may be extended as above to $m^{t***}:Y^{**}\times X^{**}\rightarrow Z^{**}$.
 The mapping $m^{t***t}:X^{**}\times Y^{**}\rightarrow Z^{**}$ in general is not equal to $m^{***}$, see [1], if $m^{***}=m^{t***t}$, then $m$ is called Arens regular. The mapping $y^{\prime\prime}\rightarrow m^{t***t}(x^{\prime\prime},y^{\prime\prime})$ is $weak^*-to-weak^*$ continuous for every $y^{\prime\prime}\in Y^{**}$, but the mapping $x^{\prime\prime}\rightarrow m^{t***t}(x^{\prime\prime},y^{\prime\prime})$ from $X^{**}$ into $Z^{**}$ is not in general  $weak^*-to-weak^*$ continuous for every $y^{\prime\prime}\in Y^{**}$. So we define the second topological center of $m$ as
$$Z_2(m)=\{y^{\prime\prime}\in Y^{**}:~~x^{\prime\prime}\rightarrow m^{t***t}(x^{\prime\prime},y^{\prime\prime})~~is~~weak^*-to-weak^*-continuous\}.$$
It is clear that $m$ is Arens regular if and only if $Z_1(m)=X^{**}$ or $Z_2(m)=Y^{**}$. Arens regularity of $m$ is equivalent to the following
$$\lim_i\lim_j<z^\prime,m(x_i,y_j)>=\lim_j\lim_i<z^\prime,m(x_i,y_j)>,$$
whenever both limits exist for all bounded sequences $(x_i)_i\subseteq X$ , $(y_i)_i\subseteq Y$ and $z^\prime\in Z^*$, see [5, 20].\\
 The regularity of a normed algebra $A$ is defined to be the regularity of its algebra multiplication when considered as a bilinear mapping. Let $a^{\prime\prime}$ and $b^{\prime\prime}$ be elements of $A^{**}$, the second dual of $A$. By $Goldstin^,s$ Theorem [4, P.424-425], there are nets $(a_{\alpha})_{\alpha}$ and $(b_{\beta})_{\beta}$ in $A$ such that $a^{\prime\prime}=weak^*-\lim_{\alpha}a_{\alpha}$ ~and~  $b^{\prime\prime}=weak^*-\lim_{\beta}b_{\beta}$. So it is easy to see that for all $a^\prime\in A^*$,
$$\lim_{\alpha}\lim_{\beta}<a^\prime,m(a_{\alpha},b_{\beta})>=<a^{\prime\prime}b^{\prime\prime},a^\prime>$$ and
$$\lim_{\beta}\lim_{\alpha}<a^\prime,m(a_{\alpha},b_{\beta})>=<a^{\prime\prime}ob^{\prime\prime},a^\prime>,$$
where $a^{\prime\prime}.b^{\prime\prime}$ and $a^{\prime\prime}ob^{\prime\prime}$ are the first and second Arens products of $A^{**}$, respectively, see [6, 17, 20].\\
The mapping $m$ is left strongly Arens irregular if $Z_1(m)=X$ and $m$ is right strongly Arens irregular if $Z_2(m)=Y$.\\
Regarding $A$ as a Banach $A-bimodule$, the operation $\pi:A\times A\rightarrow A$ extends to $\pi^{***}$ and $\pi^{t***t}$ defined on $A^{**}\times A^{**}$. These extensions are known, respectively, as the first (left) and the second (right) Arens products, and with each of them, the second dual space $A^{**}$ becomes a Banach algebra. In this situation, we shall also simplify our notations. So the first (left) Arens product of $a^{\prime\prime},b^{\prime\prime}\in A^{**}$ shall be simply indicated by $a^{\prime\prime}b^{\prime\prime}$ and defined by the three steps:
 $$<a^\prime a,b>=<a^\prime ,ab>,$$
  $$<a^{\prime\prime} a^\prime,a>=<a^{\prime\prime}, a^\prime a>,$$
  $$<a^{\prime\prime}b^{\prime\prime},a^\prime>=<a^{\prime\prime},b^{\prime\prime}a^\prime>.$$
 for every $a,b\in A$ and $a^\prime\in A^*$. Similarly, the second (right) Arens product of $a^{\prime\prime},b^{\prime\prime}\in A^{**}$ shall be  indicated by $a^{\prime\prime}ob^{\prime\prime}$ and defined by :
 $$<a oa^\prime ,b>=<a^\prime ,ba>,$$
  $$<a^\prime oa^{\prime\prime} ,a>=<a^{\prime\prime},a oa^\prime >,$$
  $$<a^{\prime\prime}ob^{\prime\prime},a^\prime>=<b^{\prime\prime},a^\prime ob^{\prime\prime}>.$$
  for all $a,b\in A$ and $a^\prime\in A^*$.\\
  The regularity of a normed algebra $A$ is defined to be the regularity of its algebra multiplication when considered as a bilinear mapping. Let $a^{\prime\prime}$ and $b^{\prime\prime}$ be elements of $A^{**}$, the second dual of $A$. By $Goldstine^,s$ Theorem [4, P.424-425], there are nets $(a_{\alpha})_{\alpha}$ and $(b_{\beta})_{\beta}$ in $A$ such that $a^{\prime\prime}=weak^*-\lim_{\alpha}a_{\alpha}$ ~and~  $b^{\prime\prime}=weak^*-\lim_{\beta}b_{\beta}$. So it is easy to see that for all $a^\prime\in A^*$,
$$\lim_{\alpha}\lim_{\beta}<a^\prime,\pi(a_{\alpha},b_{\beta})>=<a^{\prime\prime}b^{\prime\prime},a^\prime>$$ and
$$\lim_{\beta}\lim_{\alpha}<a^\prime,\pi(a_{\alpha},b_{\beta})>=<a^{\prime\prime}ob^{\prime\prime},a^\prime>,$$
where $a^{\prime\prime}b^{\prime\prime}$ and $a^{\prime\prime}ob^{\prime\prime}$ are the first and second Arens products of $A^{**}$, respectively, see [6, 17, 20].\\
We find the usual first and second topological center of $A^{**}$, which are
  $$Z_1(A^{**})=Z^\ell_1(A^{**})=\{a^{\prime\prime}\in A^{**}: b^{\prime\prime}\rightarrow a^{\prime\prime}b^{\prime\prime}~ is~weak^*-to-weak^*~continuous\},$$
   $$Z_2(A^{**})=Z_2^r(A^{**})=\{a^{\prime\prime}\in A^{**}: a^{\prime\prime}\rightarrow a^{\prime\prime}ob^{\prime\prime}~ is~weak^*-to-weak^*~continuous\}.$$\\
 An element $e^{\prime\prime}$ of $A^{**}$ is said to be a mixed unit if $e^{\prime\prime}$ is a
right unit for the first Arens multiplication and a left unit for
the second Arens multiplication. That is, $e^{\prime\prime}$ is a mixed unit if
and only if,
for each $a^{\prime\prime}\in A^{**}$, $a^{\prime\prime}e^{\prime\prime}=e^{\prime\prime}o a^{\prime\prime}=a^{\prime\prime}$. By
[4, p.146], an element $e^{\prime\prime}$ of $A^{**}$  is  mixed
      unit if and only if it is a $weak^*$ cluster point of some BAI $(e_\alpha)_{\alpha \in I}$  in
      $A$.\\
Let now $B$ be a Banach $A-bimodule$, and let
$$\pi_\ell:~A\times B\rightarrow B~~~and~~~\pi_r:~B\times A\rightarrow B.$$
be the right and left module actions of $A$ on $B$. Then $B^{**}$ is a Banach $A^{**}-bimodule$ with module actions
$$\pi_\ell^{***}:~A^{**}\times B^{**}\rightarrow B^{**}~~~and~~~\pi_r^{***}:~B^{**}\times A^{**}\rightarrow B^{**}.$$
Similarly, $B^{**}$ is a Banach $A^{**}-bimodule$ with module actions\\
$$\pi_\ell^{t***t}:~A^{**}\times B^{**}\rightarrow B^{**}~~~and~~~\pi_r^{t***t}:~B^{**}\times A^{**}\rightarrow B^{**}.$$\\

\begin{center}
\section{ \bf Cohomological groups of Banach algebras  }
\end{center}
Let $B$ be a Banach $A-bimodule$ and let $n\geq 1$. Suppose that $B^{(n)}$ is an  $n-th~dual$ of  $B$. Then $B^{(n)}$ is also Banach $A-bimodule$, that is, for every $a\in A$, $b^{(n)}\in B^{(n)}$ and $b^{(n-1)}\in B^{(n-1)}$, we define
 $$<b^{(n)}a,b^{(n-1)}>= <b^{(n)},ab^{(n-1)}>,$$
 and
 $$<ab^{(n)},b^{(n-1)}>= <b^{(n)},b^{(n-1)}a>.$$
 In the following theorem, we extend  Theorem 1.9 from [6] into general situation.\\\\

\noindent{\it{\bf Theorem 2-1.}} Let $B$ be a Banach $A-bimodule$ and let $n\geq 1$. If $H^1(A,B^{(n+2)})=0$, then $H^1(A,B^{(n)})=0$.\\
\begin{proof} Let $D\in Z^1(A,B^{(n)})$ and suppose that  $i: B^{(n)}\rightarrow B^{(n+2)}$ is the canonical linear mapping as $A-bimodule$ homomorphism. Take $\widetilde{D}=ioD$. Then we can be viewed $\widetilde{D}$ as an element of
$Z^1(A,B^{(n+2)})$. Since $H^1(A,B^{(n+2)})=0$, there exist $b^{(n+2)}\in B^{(n+2)}$ such that for every $a\in A$, we have $$\widetilde{D}(a)=ab^{(n+2)}-b^{(n+2)}a.$$
Set a $A-linear$ mapping $P$ from $B^{(n+2)}$ into $B^{(n)}$ such that $Poi=I_{B^{(n)}}$. Then we have $Po\widetilde{D}=(Poi)oD=D$, and so for every $a\in A$, we conclude that $D(a)=Po\widetilde{D}(a)=aP(b^{(n+2)})-P(b^{(n+2)})a$. It follows that $D\in N^1(A,B^{(n)})$. Consequently we have $H^1(A,B^{(n)})=0$.\\
\end{proof}

 Suppose that $A$ is a Banach algebra and $B$ is a Banach $A-bimodule$. According to [6, pp.27 and 28], $B^{**}$ is a Banach $A^{**}-bimodule$, where  $A^{**}$ is equipped with the first Arens product. So we can define the topological centers of module actions.\\
Thus,  for a Banach $A-bimodule$ $B$,  we  define the topological centers of the  left and right module actions of $A^{**}$ on $B^{**}$ as follows:

$${Z}^\ell_{A^{**}}(B^{**})={Z}(\pi_r)=\{b^{\prime\prime}\in B^{**}:~the~map~~a^{\prime\prime}\rightarrow \pi_r^{***}(b^{\prime\prime}, a^{\prime\prime})~:~A^{**}\rightarrow B^{**}$$$$~is~~~weak^*-to-weak^*~continuous\}$$
$${Z}^\ell_{B^{**}}(A^{**})={Z}(\pi_\ell)=\{a^{\prime\prime}\in A^{**}:~the~map~~b^{\prime\prime}\rightarrow \pi_\ell^{***}(a^{\prime\prime}, b^{\prime\prime})~:~B^{**}\rightarrow B^{**}$$$$~is~~~weak^*-to-weak^*~continuous\}$$
$${Z}^r_{A^{**}}(B^{**})={Z}(\pi_\ell^t)=\{b^{\prime\prime}\in B^{**}:~the~map~~a^{\prime\prime}\rightarrow \pi_\ell^{t***}(b^{\prime\prime}, a^{\prime\prime})~:~A^{**}\rightarrow B^{**}$$$$~is~~~weak^*-to-weak^*~continuous\}$$
$${Z}^r_{B^{**}}(A^{**})={Z}(\pi_r^t)=\{a^{\prime\prime}\in A^{**}:~the~map~~b^{\prime\prime}\rightarrow \pi_r^{t***}(a^{\prime\prime}, b^{\prime\prime})~:~B^{**}\rightarrow B^{**}$$$$~is~~~weak^*-to-weak^*~continuous\}.$$\\
In every parts of this paper, for left or right Banach $A-module$ $B$, we take $\pi_\ell(a,b)=ab$ and  $\pi_r(b,a)=ba$, for each $a\in A$ and $b\in B$.\\
Let $A^{(n)}$ and  $B^{(n)}$  be $n-th~dual$ of $A$ and $B$, respectively. By [25, page 4132-4134], if $n\geq 0$ is an even number, then  $B^{(n)}$ is a Banach $A^{(n)}-bimodule$. Then for $n\geq 2$,   we define  $B^{(n)}B^{(n-1)}$ as a subspace of $A^{(n-1)}$, that is, for all $b^{(n)}\in B^{(n)}$,  $b^{(n-1)}\in B^{(n-1)}$ and  $a^{(n-2)}\in A^{(n-2)}$ we define
$$<b^{(n)}b^{(n-1)},a^{(n-2)}>=<b^{(n)},b^{(n-1)}a^{(n-2)}>.$$
If $n$ is odd number, then for $n\geq 1$,   we define  $B^{(n)}B^{(n-1)}$ as a subspace of $A^{(n)}$, that is, for all $b^{(n)}\in B^{(n)}$,  $b^{(n-1)}\in B^{(n-1)}$ and  $a^{(n-1)}\in A^{(n-1)}$ we define
$$<b^{(n)}b^{(n-1)},a^{(n-1)}>=<b^{(n)},b^{(n-1)}a^{(n-1)}>.$$
and if $n=0$, we take $A^{(0)}=A$ and $B^{(0)}=B$.\\
We also define the topological centers of module actions of $A^{(n)}$ on  $B^{(n)}$ as follows

$${Z}^\ell_{A^{(n)}}(B^{(n)})=\{b^{(n)}\in B^{(n)}:~the~map~~a^{(n)}\rightarrow b^{(n)} a^{(n)}~:~A^{(n)}\rightarrow B^{(n)}$$$$~is~~~weak^*-to-weak^*~continuous\}$$
$${Z}^\ell_{B^{(n)}}(A^{(n)})=\{a^{(n)}\in A^{(n)}:~the~map~~b^{(n)}\rightarrow a^{(n)} b^{(n)}~:~B^{(n)}\rightarrow B^{(n)}$$$$~is~~~weak^*-to-weak^*~continuous\}.$$\\

In the following, for every $n\geq 1$, using of the Arens regularity of left module action of $A^{(2n)}$ on $B^{(2n)}$, we extend every derivation $D:A\rightarrow B^{(2n)}$ into derivation $\widetilde{D}$ from $A^{(2n)}$ into $B^{(2n)}$ such that $\widetilde{D}(a)=D(a)$ for every $a\in A$ and we have some conclusion in cohomological groups.\\\\

\noindent{\it{\bf Theorem 2-2.}} Let $B$ be a Banach $A-bimodule$ and  $D:A\rightarrow B^{(2n)}$ be a continuous derivation. Assume that $Z^\ell_{A^{(2n)}}(B^{(2n)})=B^{(2n)}$. Then there is a continuous derivation $\widetilde{D}: A^{(2n)}\rightarrow B^{(2n)}$ such that $\widetilde{D}(a)=D(a)$ for every $a\in A$.
\begin{proof} By using [6, Proposition 1.7], the linear mapping $D^{\prime\prime}:A^{**}\rightarrow B^{(2n+2)}$ is a continuous derivation. Take $X=B^{(2n-2)}$. Since $Z_{A^{(2n)}}(X^{**})=Z_{A^{(2n)}}(B^{(2n)})=B^{(2n)}=X^{**}$, by [6, Proposition 1.8], the canonical projection $P:  X^{(4)}\rightarrow X^{**}$ is a $A^{**}-bimodule $ morphism. Set $\widetilde{D}=PoD^{\prime\prime}$. Then $\widetilde{D}$ is a continuous derivation from $A^{**}$ into $B^{(2n)}$. Now by replacing $A^{**}$ by $A$ and repeating of the proof, we obtain the result. \\ \end{proof}

\noindent{\it{\bf Corollary 2-3.}} Let $B$ be a Banach $A-bimodule$ and $n\geq 0$. If $Z^\ell_{A^{(2n)}}(B^{(2n)})=B^{(2n)}$ and $H^1(A^{(2n+2)},B^{(2n+2)})=0$, then  $H^1(A,B^{(2n)})=0$.
\begin{proof} By using [6, Proposition 1.7] and proceeding theorem the result is hold.\\ \end{proof}

Let $A$ be a Banach algebra and $n\geq 0$. Then $A$ is $n-weakly$ amenable if $H^1(A,A^{(n)})=0$. $A$ is permanently weakly amenable if $A$ is $n-weakly$ amenable for each $n\geq 0$.\\\\

\noindent{\it{\bf Corollary 2-4}} [6]. Let $A$ be a Banach algebra such that $A^{(2n)}$ is Arens regular and
 $H^1(A^{(2n+2)}),A^{(2n+2)})=0$ for each $n\geq 0$. Then $A$ is $2n-weakly$ amenable.\\\\

Assume that $A$ is Banach algebra and $n\geq 0$. We define $A^{[n]}$ as a subset of $A$ as follows
$$A^{[n]}=\{a_1a_2...a_n:~a_1,a_2,...a_n\in A\}.$$
We write $A^n$  the linear span of $A^{[n]}$ as a subalgebra of $A$.\\\\

\noindent{\it{\bf Theorem 2-5.}} Let $A$ be a Banach algebra and $n\geq 0$. Let $A^{[2n]}$ dense in $A$ and suppose that $B$ is a Banach $A-bimodule$. Assume that $AB^{**}$ and $B^{**}A$ are subsets of $B$. If $H^1(A,B^*)=0$, then  $H^1(A ,B^{(2n+1)})=0$.
\begin{proof} For $n=0$ the result is clear.  Let $B^\perp$ be the space of functionals in $A^{(2n+1)}$ which annihilate $i(B)$ where $i:B\rightarrow
B^{(2n)}$ is a natural canonical mapping. Then, by using [26, lemma 1], we have the following equality.
$$B^{(2n+1)}=i(B^*)\oplus B^\bot ,$$
as Banach $A-bimodules$, and so
$$H^1(A,B^{(2n+1)})=H^1(A,i(B^*))\oplus H^1(A,B^\bot ) .$$ Without lose generality, we replace $i(B^*)$ by $B^*$.
Since $H^1(A,B^*)=0$, it is enough to show that $H^1(A,B^\bot )=0$.\\
Now, take the linear mappings $L_a$ and $R_a$ from $B$ into itself by $L_a(b)=ab$ and $R_a(b)=ba$ for every $a\in A$.
Since $AB^{**}\subseteq B$ and $B^{**}A\subseteq B$,  $L^{**}_a(b^{\prime\prime})=ab^{\prime\prime}$ and $R^{**}_a(b^{\prime\prime})=b^{\prime\prime}a$ for every $a\in A$, respectively. Consequently, $L_a$ and $R_a$ from $B$ into itself are weakly compact. It follows that for each $a\in A$ the linear mappings $L^{(2n)}_a$ and   $R^{(2n)}_a$ from $B^{(n)}$ into $B^{(n)}$ are weakly compact and for every $b^{(2n)}\in B^{(2n)}$, we have $L^{(2n)}_a(b^{(2n)})=ab^{(2n)}\in B^{(2n-2)}$ and $R^{(2n)}_a(b^{(2n)})=b^{(2n)}a\in B^{(2n-2)}$. Set $a_1, a_2, ...,a_n\in A$ and
$b^{(2n)}\in B^{(2n)}$. Then $a_1 a_2 ...a_nb^{(2n)}$ and $b^{(2n)}a_1 a_2 ...a_n$ are belong to $B$.
Suppose that $D\in Z^1(A,B^\bot )$ and let $a, x\in A^{[n]}$. Then for every $b^{(2n)}\in B^{(2n)}$, since $xb^{(2n)}, b^{(2n)}a \in B$, we have the following equality
$$<D(ax), b^{(2n)}>=<aD(x), b^{(2n)}>+<D(a)x, b^{(2n)}>$$$$=<D(x), b^{(2n)}a>+<D(a), xb^{(2n)}>=0.$$
It follows that $D\mid_{A^{[2n]}}$. Since $A^{[2n]}$ dense in $A$, $D=0$. Hence $H^1(A,B^\bot )=0$ and result follows.\\
\end{proof}

\noindent{\it{\bf Corollary 2-6.}}\\ i) Let $A$ be a Banach algebra with bounded left approximate identity [=$LBAI$], and let  $B$ be a Banach $A-bimodule$.  Suppose that $AB^{**}$ and $B^{**}A$ are subset of $B$. Then if $H^1(A,B^*)=0$, it follows that  $H^1(A,B^{(2n+1)})=0$.\\
ii) Let $A$ be an amenable Banach algebra and $B$ be a Banach $A-bimodule$. If $AB^{**}$ and $B^{**}A$ are subset of $B$, then $H^1(A,B^{(2n+1)})=0$.\\\\

\noindent{\it{\bf Example 2-7.}} Assume that $G$ is a compact group. Then \\
i) we know that $L^1(G)$ is $M(G)-bimodule$  and $L^1(G)$ is an ideal in the second dual of $M(G)$, $M(G)^{**}$. By using [18, corollary 1.2], we have $H^1(L^1(G),M(G))=0$. Then  for every $n\geq 1$, by using proceeding corollary, we conclude that
$$H^1(L^1(G),M(G)^{(2n+1)})=0.$$
ii) we have $L^1(G)$ is an ideal in its second dual , $L^1(G)^{**}$. By using [15], we know that $L^1(G)$ is a weakly amenable. Then by proceeding corollary, $L^1(G)$ is $(2n+1)-weakly$  amenable.\\\\

\noindent{\it{\bf Corollary 2-8.}} Let $A$ be a Banach algebra and let $A^{[2n]}$ be dense in $A$. Suppose that $AB^{**}$ and $B^{**}A$ are subset of $B$. Then the following are equivalent.
\begin{enumerate}
\item ~$H^1(A,B^*)=0$.
\item ~$H^1(A,B^{(2n+1)})=0$ for some $n\geq 0$.

\item ~~$H^1(A,B^{(2n+1)})=0$ for each $n\geq 0$.\\
\end{enumerate}

\noindent{\it{\bf Corollary 2-9}} [6]. Let $A$ be a weakly amenable Banach algebra such that $A$ is an ideal in $A^{**}$. Then $A$ is $(2n+1)-weakly$ amenable for each $n\geq 0$.
\begin{proof} By using Proposition 1.3 from [6] and proceeding theorem, result is hold.\\ \end{proof}

Let $B$ be a dual Banach algebra and $dimB<\infty$. Then by using Proposition 2.6.24 from [5], we know that $\mathcal{N}(B)$, the collection of all operators from $B$ into $B$, is an ideal in $\mathcal{N}(B)^{**}$. By using Corollary 4.3.6 from [21], $N(B)$ is amenable, and so it also is weakly amenable. Consequently, by using the proceeding corollary, $\mathcal{N}(B)$ is $(2n+1)-weakly$ amenable for every $n\geq 0$.\\
We know that every von Neumann algebra $A$ is an weakly amenable Banach algebra, see [21]. Now if $A$ is an ideal in its second dual, $A^{**}$, then by using proceeding corollary, $A$ is $(2n+1)-weakly$ amenable Banach algebra for each $n\geq 0$.\\

Assume that $A$ and $B$ are Banach algebra. Then $A\oplus B$ , with norm $$\parallel(a,b)\parallel=\parallel a\parallel+\parallel b\parallel,$$ and product $(a_1,b_1)(a_2,b_2)=(a_1a_2, b_1b_2)$ is a Banach algebra. It is clear that if $X$ is a Banach $A~and ~B-bimodule$, then $X$ is a Banach $A\oplus B-bimodule$.\\
In the following, we investigated the relationships between the cohomological group of  $A\oplus B$ and cohomological groups of $A$ and $B$.\\\\

\noindent{\it{\bf Theorem 2-10.}} Suppose that $A$ and $B$ are Banach algebras. Let $X$ be a Banach $A~and ~B-bimodule$. Then,   $H^1(A\oplus B, X)=0$ if and only if $H^1(A,X)=H^1(B,X)=0$.
\begin{proof} Suppose that $H^1(A\oplus B, X)=0$. Assume that $D_1\in Z^1(A,X)$ and $D_2\in Z^1(B,X)$. Take $D=(D_1,D_2)$. Then for every $a_1,a_2\in A$ and $b_1, b_2\in B$, we have
$$D((a_1,b_1)(a_2,b_2))=D(a_1a_2,b_1b_2)=(D_1(a_1a_2),D_2(b_1b_2))$$
$$=(a_1D_1(a_2)+D_1(a_1)a_2,b_1D_2(b_2)+D_2(b_1)b_2)$$
$$=(a_1D_1(a_2),b_1D_2(b_2))+(D_1(a_1)a_2+D_2(b_1)b_2)$$
$$=(a_1,b_1)(D_1(a_2),D_2(b_2))+(D_1(a_1),D_2(b_1))(a_2,b_2)$$
$$=(a_1,b_1)D(a_2,b_2)+D(a_1,b_1)(a_2,b_2).$$
It follows that $D\in Z^1(A\oplus B,X)$.
Since $H^1(A\oplus B, X)=0$, there is $x\in X$ such that $D=\delta_{x}$  where  $\delta_{x}\in N^1(A\oplus B,X)$. Since
$\delta_{x}=(\delta^1_{x},\delta^2_{x})$ where $\delta^1_{x}\in N^1(A,X)$ and $\delta^2_{x}\in N^1(B,X)$,  we have $D_1=\delta^1_{x}$ and $D_2=\delta^2_{x}$. Thus we have
$H^1(A,X)=H^1(B,X)=0$.\\
For the converse, take $A$ as an ideal in $A\oplus B$, and so by using Proposition 2.8.66 from [5], proof hold.\\ \end{proof}

Let $G$ be a locally compact group and $X$ be a Banach $L^1(G)-bimodule$. Then by [6, pp.27 and 28], $X^{**}$ is a Banach $L^1(G)^{**}-bimodule$.  Since $L^1(G)^{**}=LUC(G)^*\oplus LUC(G)^\bot$, by using proceeding theorem, we have  $H^1(L^1(G)^{**},X^{**})=0$ if and only if $H^1(LUC(G)^*,X^{**})=H^1(LUC(G)^\bot,X^{**})=0$.\\
On the other hand, we know that $L^1(G)^{**}=L^1(G)\oplus C_0(G)^\perp$. By [15], we know that, $H^1(L^1(G), L^\infty(G))=0$. By using proceeding theorem,  $$H^1(L^1(G)^{**}, L^\infty(G))=0,$$ if and only if  $H^1(C_0(G)^\perp, L^\infty(G))=0$.\\\\

\noindent{\it{\bf Corollary 2-11.}} \\
i) Suppose that $A$ and $B$ are Banach algebras. Let $X$ be a Banach $A~and ~B-bimodule$. Then   $A\oplus B$ is an amenable Banach algebra if and only if $A$ and $B$ are amenable Banach algebras.\\
ii) Let $A$ be a Banach algebra and $n\geq 1$. Then  $H^1(\oplus_{i=1}^nA,A^*)=0$ if and only if $A$ is weakly amenable.\\\\

Assume that $B$ is a Banach $A-bimodule$. In the following, we will study the relationships between two cohomological groups $H^1(A,B^*)$ and  $H^1(A^{**},B^{***})$. If  $H^1(A,B^*)=0$, we want to know that dose $H^1(A^{**},B^{***})=0$ and when its converse hold? If we give an answer to these questions, then we want to establish the results of them in the weak amenability of Banach algebras. As some applications of these discussion in the group algebras, for a compact group $G$, we show that every $weak^*-to-weak$ continuous derivation from $L^1(G)^{**}$ into $M(G)$ is inner. \\\\

\noindent{\it{\bf Theorem 2-12.}} Let $B$ be a Banach $A-bimodule$ and suppose that every derivation from $A$ into $B^*$ is weakly compact. If $Z^\ell_{A^{**}}(B^{**})=B^{**}$ and $H^1(A^{**},B^{***})=0$, then $H^1(A,B^*)=0$.
\begin{proof} Suppose that $D\in Z^1(A,B^*)$. Since $D:A\rightarrow B^*$ is weakly compact $D^{\prime\prime}(A^{**})\subseteq B^*$. Assume that $a^{\prime\prime}, x^{\prime\prime}\in A^{**}$ and $(a_\alpha)_\alpha, (x_\beta)_\beta\subseteq A$ such that $a_\alpha^{} \stackrel{w^*} {\rightarrow} a^{\prime\prime}$ and $x_\beta \stackrel{w^*} {\rightarrow} x^{\prime\prime}$ in $A^{**}$. Then, since $Z^\ell_{A^{**}}(B^{**})=B^{**}$, for every $b^{\prime\prime}\in B^{**}$, we have
$$\lim_\alpha <a_\alpha D^{\prime\prime}(x^{\prime\prime}),b^{\prime\prime}>=
\lim_\alpha < D^{\prime\prime}(x^{\prime\prime}),b^{\prime\prime}a_\alpha>=
<b^{\prime\prime}a^{\prime\prime}, D^{\prime\prime}(x^{\prime\prime})>$$$$=
 <b^{\prime\prime},a^{\prime\prime} D^{\prime\prime}(x^{\prime\prime})>=
 <a^{\prime\prime} D^{\prime\prime}(x^{\prime\prime}),b^{\prime\prime}>.$$
Consequently, we have
$$D^{\prime\prime}(a^{\prime\prime}x^{\prime\prime})=\lim_\alpha\lim_\beta D(a_\alpha x_\beta)=
\lim_\alpha\lim_\beta (D(a_\alpha) x_\beta +a_\alpha D(x_\beta))$$$$=D^{\prime\prime}(a^{\prime\prime})x^{\prime\prime}
+a^{\prime\prime}D^{\prime\prime}(x^{\prime\prime}).$$
It follows that $D^{\prime\prime}\in Z^1(A^{**},B^{***})$. Since $H^1(A^{**},B^{***})=0$, there is $b^{\prime\prime\prime}\in  B^{***}$ such that $D^{\prime\prime}=\delta_{b^{\prime\prime\prime}}$. Then we have $D=D^{\prime\prime}\mid_A=\delta_{b^{\prime\prime\prime}}\mid_A$. So for every $a\in A$, we have $D(a)=\delta_{b^{\prime\prime\prime}}(a)=ab^{\prime\prime\prime}-b^{\prime\prime\prime}a$. If we take $b^\prime=b^{\prime\prime\prime}\mid_{B}$, then $D(a)=ab^\prime-b^\prime a=\delta_{b^\prime}(a)$. We conclude that $H^1(A,B^*)=0$.\\ \end{proof}

Let $A$ be a Banach algebra, and let $B$ be a Banach $A-bimodule$. A continuous derivation $D:A\rightarrow B$ is approximately inner [resp. weakly approximately inner], if there exists a bounded net $(b_\alpha)_\alpha \subseteq B$ such that $D(a)=\lim_\alpha (ab_\alpha-b_\alpha a)$  in $B$ [resp. $D(a)=w-\lim_\alpha (ab_\alpha-b_\alpha a)$ in $B$].\\
Assume that $(b^\prime_\alpha)_\alpha \subseteq B^*$ such that $b^\prime_\alpha \stackrel{w^*} {\rightarrow} b^{\prime\prime\prime}$ in $B^{***}$. Then by assumptions of proceeding theorem, for every $a\in A$ and for every derivation from $A$ into $B^*$,   we have
$D(a)=D^{\prime\prime}\mid_A(a)=\delta_{b^{\prime\prime\prime}}\mid_A(a)=ab^{\prime\prime\prime}-b^{\prime\prime\prime
}a=w^*-\lim_\alpha (ab^\prime_\alpha-b^\prime_\alpha a)$ in $B^{***}$. Since $D(a)\in B^*$, $D(a)=w-\lim_\alpha (ab^\prime_\alpha-b^\prime_\alpha a)$ in $B^*$, and so $D$ is weakly approximately inner in $B$.\\\\

\noindent{\it{\bf Theorem 2-13.}}  Let $B$ be a Banach $A-bimodule$ and suppose that $D^{\prime\prime}(A^{**})\subseteq B^*$. If $Z^\ell_{A^{**}}(B^{**})=B^{**}$ and $H^1(A^{**},B^{***})=0$, then $H^1(A,B^*)=0$.
\begin{proof} Proof is similar to proceeding theorem.\\ \end{proof}

A functional $a^\prime$ in $A^*$ is said to be $wap$ (weakly almost
 periodic) on $A$ if the mapping $a\rightarrow a^\prime a$ from $A$ into
 $A^{*}$ is weakly compact.  In  [20], Pym showed that  this definition to the equivalent following condition\\
 For any two net $(a_{\alpha})_{\alpha}$ and $(b_{\beta})_{\beta}$
 in $\{a\in A:~\parallel a\parallel\leq 1\}$, we have\\
$$\\lim_{\alpha}\\lim_{\beta}<a^\prime,a_{\alpha}b_{\beta}>=\\lim_{\beta}\\lim_{\alpha}<a^\prime,a_{\alpha}b_{\beta}>,$$
whenever both iterated limits exist. The collection of all $wap$
functionals on $A$ is denoted by $wap(A)$. Also we have
$a^{\prime}\in wap(A)$ if and only if $<a^{\prime\prime}b^{\prime\prime},a^\prime>=<a^{\prime\prime}ob^{\prime\prime},a^\prime>$ for every $a^{\prime\prime},~b^{\prime\prime} \in
A^{**}$. \\
Let $B$ be a Banach left $A-module$. Then, $b^\prime\in B^*$ is said to be left weakly almost periodic functional if the set $\{\pi^*_\ell(b^\prime,a):~a\in A,~\parallel a\parallel\leq 1\}$ is relatively weakly compact. We denote by $wap_\ell(B)$ the closed subspace of $B^*$ consisting of all the left weakly almost periodic functionals in $B^*$.\\
The definition of the right weakly almost periodic functional ($=wap_r(B)$) is the same.\\
By  [20], $b^\prime\in wap_\ell(B)$ is equivalent to the following $$<\pi_\ell^{***}(a^{\prime\prime},b^{\prime\prime}),b^\prime>=
<\pi_\ell^{t***t}(a^{\prime\prime},b^{\prime\prime}),b^\prime>$$
for all $a^{\prime\prime}\in A^{**}$ and $b^{\prime\prime}\in B^{**}$.
Thus, we can write \\
$$wap_\ell(B)=\{ b^\prime\in B^*:~<\pi_\ell^{***}(a^{\prime\prime},b^{\prime\prime}),b^\prime>=
<\pi_\ell^{t***t}(a^{\prime\prime},b^{\prime\prime}),b^\prime>~~$$$$for~~all~~a^{\prime\prime}\in A^{**},~b^{\prime\prime}\in B^{**}\}.$$\\\\

\noindent{\it{\bf Corollary 2-14.}}  Let $B$ be a Banach $A-bimodule$ and let every derivation\\ $D:A\rightarrow B^*$,
satisfies $D^{\prime\prime}(A^{**})\subseteq wap_\ell(B)$. If $H^1(A^{**},B^{***})=0$, then\\ $H^1(A,B^*)=0$.\\

In Corollary 2-14, if we take $B=A$, then we obtain Theorem 2.1 from [10].\\

\noindent{\it{\bf Corollary 2-15.}}  Let $B$ be a Banach $A-bimodule$ and suppose that for every derivation $D:A\rightarrow B^*$, we have $D^{\prime\prime}(A^{**})B^{**}\subseteq A^*$. If $H^1(A^{**},B^{***})=0$, then $$H^1(A,B^*)=0.$$\\

\noindent{\it{\bf Theorem 2-16.}}  Let $B$ be a Banach $A-bimodule$ and suppose that $AA^{**}\subseteq A$ and
$B^{**}A=B^{**}$. If $H^1(A^{**},B^{***})=0$, then $H^1(A,B^*)=0$.
\begin{proof} Let   $(a_{\alpha}^{\prime\prime})_{\alpha}\subseteq A^{**}$ such that  $a^{\prime\prime}_{\alpha} \stackrel{w^*} {\rightarrow}a^{\prime\prime}$ in $A^{**}$. Assume that $b^{\prime\prime}\in B^{**}$. Since $B^{**}A=B^{**}$, there are $a\in A$ and $y^{\prime\prime}\in B^{**}$ such that $b^{\prime\prime}=y^{\prime\prime}a$. We know that $a\in A\subseteq Z^\ell_{B^{**}}(A^{**})$, and so $aa^{\prime\prime}_{\alpha} \stackrel{w^*} {\rightarrow}aa^{\prime\prime}$ in $A^{**}$. Since $(aa^{\prime\prime}_{\alpha})_{\alpha}\subseteq A$ and $aa^{\prime\prime}\in A$, we have $aa^{\prime\prime}_{\alpha} \stackrel{w} {\rightarrow}aa^{\prime\prime}$ in $A$.
Then for every $a^{\prime\prime}\in A^{**}$, we have the following equalities
$$<D^{\prime\prime}(x^{\prime\prime})b^{\prime\prime},a^{\prime\prime}_{\alpha}>=
<D^{\prime\prime}(x^{\prime\prime})y^{\prime\prime}a,a^{\prime\prime}_{\alpha}>=
<D^{\prime\prime}(x^{\prime\prime})y^{\prime\prime},aa^{\prime\prime}_{\alpha}>\rightarrow
<D^{\prime\prime}(x^{\prime\prime})y^{\prime\prime},aa^{\prime\prime}>$$$$=
<D^{\prime\prime}(x^{\prime\prime})b^{\prime\prime},a^{\prime\prime}>.$$
It follows that $D^{\prime\prime}(x^{\prime\prime})b^{\prime\prime}\in (A^{**},weak^*)^*=A^*$. So, by using Corollary 2-15, proof is hold.\\ \end{proof}

\noindent{\it{\bf Theorem 2-17.}} Let $B$ be a Banach $A-bimodule$  and let every derivation $D:A^{**}\rightarrow B^*$ is $weak^*-to-weak^*$ continuous. Then if $Z^\ell_{B^{**}}(A^{**})=A^{**}$ and $H^1(A,B^*)=0$, it follows that  $H^1(A^{**},B^*)=0$.
\begin{proof} Let $D:A^{**}\rightarrow B^*$ be a derivation. Then $D\mid_A:A\rightarrow B^*$ is a derivation. Since $H^1(A,B^*)=0$, there is $b^\prime\in B^*$ such that $D\mid_A=\delta_{b^\prime}$. Suppose that $a^{\prime\prime}\in A^{**}$ and
$(a_\alpha^{})_\alpha\subseteq A^{}$ such that $a_\alpha^{} \stackrel{w^*} {\rightarrow} a^{\prime\prime}$ in $A^{**}$. Then
$$D(a^{\prime\prime})=w^*-\lim_\alpha D\mid_A(a_\alpha)=w^*-\lim_\alpha\delta_{b^\prime}(a_\alpha)=w^*-\lim_\alpha(a_\alpha {b^\prime}-{b^\prime}a_\alpha)$$$$=a^{\prime\prime}{b^\prime}-{b^\prime}a^{\prime\prime}.$$
Now, we show that $b^\prime a^{\prime\prime}\in B^*$. Assume that  $(b^{\prime\prime}_{\beta})_{\beta}\in B^{**}$ such that $b^{\prime\prime}=w^*-\lim_{\beta}b^{\prime\prime}_{\beta}$. Then since $Z^\ell_{B^{**}}(A^{**})=A^{**}$, we have
$$<b^\prime a^{\prime\prime},b^{\prime\prime}_{\beta}>=< a^{\prime\prime}b^{\prime\prime}_{\beta},b^\prime>\rightarrow
< a^{\prime\prime}b^{\prime\prime},b^\prime>=<b^\prime a^{\prime\prime},b^{\prime\prime}>.$$
Thus, we conclude that $b^\prime a^{\prime\prime}\in (B^{**},weak^*)^*=B^*$, and so
$H^1(A^{**},B^*)=0$.\\ \end{proof}

\noindent{\it{\bf Corollary 2-18.}} Let $A$ be a Arens regular  Banach algebra and let every derivation $D:A^{**}\rightarrow A^*$ is $weak^*-to-weak^*$ continuous. Then if $A$ is weakly amenable, it follows that  $H^1(A^{**},A^*)=0$.\\\\

Let $B$  be a dual Banach space, with predual  $X$ and suppose that
$$X^\perp=\{x^{\prime\prime\prime}:~x^{\prime\prime\prime}\mid_X=0~~where~~x^{\prime\prime\prime}\in X^{***}\}=\{b^{\prime\prime}:~b^{\prime\prime}\mid_X=0~~where~~b^{\prime\prime}\in B^{**}\}.$$
Then the canonical projection $P:X^{***}\rightarrow X^*$ gives a continuous linear map $P:B^{**}\rightarrow B$. Thus, we can write the following equality
$$B^{**}=X^{***}=X^*\oplus ker P=B\oplus X^\perp ,$$
as a direct sum of Banach $A-bimodule$.\\\\

\noindent{\it{\bf Theorem 2-19.}}  Let $B$ be a Banach $A-bimodule$ and  every derivation $A^{**}$ into  $ B$ is $weak^*-to-weak$ continuous. Let $A^{**}B,BA^{**}\subseteq B$. Then  the following assertions are hold.
\begin{enumerate}
\item  If $H^1(A,B)=0$, then $H^1(A^{**},B)=0$.

\item  Suppose that $A$ has a $LBAI$. Let $B$ has a predual $X$ and let $AB^*,~B^*A\subseteq X$. Then, if $H^1(A,B)=0$, it follows that  $H^1(A^{**},B^{**})=0$.

\end{enumerate}
\begin{proof}
\begin{enumerate}
\item  Proof is similar to Theorem 2-17.

\item  Set $B^{**}=B\oplus X^\perp .$ Then we have
$$H^1(A^{**},B^{**})=H^1(A^{**},B)\oplus H^1(A^{**},X^\perp ).$$
Since $H^1(A,B)=0$, by part (1), $H^1(A^{**},B)=0$. Now let $\widetilde{D}\in Z^1(A^{**},X^\perp)$ and we take $D=\widetilde{D}\mid_A$. It is clear that $D\in Z^1(A^{**},X^\perp)$. Assume that $a^{\prime\prime}, x^{\prime\prime}\in A^{**}$ and $(a_\alpha)_\alpha, (x_\beta)_\beta\subseteq A$ such that $a_\alpha^{} \stackrel{w^*} {\rightarrow} a^{\prime\prime}$ and $x_\beta \stackrel{w^*} {\rightarrow} x^{\prime\prime}$ on $A^{**}$. Since $AB^*,~B^*A\subseteq X$, for every $b^ \prime \in B^*$, by using the $weak^*-to-weak$ continuity  of $\widetilde{D}$, we have
$$<\widetilde{D}(a^{\prime\prime}x^{\prime\prime}), b^\prime>=\lim_\alpha\lim_\beta <D(a_\alpha x_\beta), b^\prime>$$$$=
\lim_\alpha\lim_\beta <(D(a_\alpha) x_\beta +a_\alpha D(x_\beta)),b^\prime>$$$$=
\lim_\alpha\lim_\beta <D(a_\alpha) x_\beta ,b^\prime>+\lim_\alpha\lim_\beta <a_\alpha D(x_\beta),b^\prime>$$
$$=
\lim_\alpha\lim_\beta <D(a_\alpha), x_\beta b^\prime>+\lim_\alpha\lim_\beta <D(x_\beta)),b^\prime a_\alpha >$$

$$=0.$$
Since $A$ has a $LBAI$, $A^{**}$ has a left unit $e^{\prime\prime}$ with respect to the first Arens product. Then  $D(x^{\prime\prime})= D(e^{\prime\prime}x^{\prime\prime})=0$, and so $D=0$.\\ \end{enumerate}
\end{proof}

Assume that $G$ is a compact group. Then we know that $L^1(G)$ is $M(G)-bimodule$  and $L^1(G)$ is an ideal in the second dual of $M(G)$, $M(G)^{**}$. By using [18, Corollary 1.2], we have $H^1(L^1(G),M(G))=0$. Then   by using proceeding corollary,  every $weak^*-to-weak$ continuous derivation from $L^1(G)^{**}$ into $M(G)$ is inner.

Consider the algebra $c_0=(c_0,.)$ is the collection of all sequences of scalars that convergence to $0$, with the some vector space operations and norm as $\ell^\infty$. We know that  $c_0$ is  a $C^*-algebra$ and Since enery $C^*-algebra$ is weakly amenable, it follows that $c_0$ is weakly amenable. Then by using proceeding theorem, every $weak^*-to-weak$ continuous derivation from $\ell_\infty$ into $\ell^1$ is inner.\\
Now let $B$ be a Banach $A-bimodule$ and for every derivation $D:A^{**}\rightarrow B$, we have $D(A^{**})\subseteq \overline{D\mid_A(A)}^w$. Assume that $D\in Z^1(A^{**},B)$. Take $\widetilde{D}=D\mid_A$. Then $\widetilde{D}\in Z^1(A,B)$. If $H^1(A,B)=0$, then there is a $b\in B$ such that $\widetilde{D}=\delta_b$.  Since $D(A^{**})\subseteq \overline{D\mid_A(A)}^w$, for every $a^{\prime\prime}\in A^{**}$, there is $(a_{\alpha})_{\alpha}\subseteq A$ such that
 $\widetilde{D} (a_{\alpha}) \stackrel{w} {\rightarrow}D(a^{\prime\prime})$. Thus we have the following equality.
 $$D(a^{\prime\prime})=w-\lim_\alpha\widetilde{D}(a_\alpha )=w-\lim_\alpha\delta_b(a_\alpha )=w-\lim_\alpha (a_\alpha b-  ba_\alpha)~in~B.$$   Thus, for every derivation $D:A^{**}\rightarrow B$, by assumption $D(A^{**})\subseteq \overline{D\mid_A(A)}^w$ and $H^1(A,B)=0$, we can write $D(a^{\prime\prime})=\lim_\alpha (a_\alpha b-  ba_\alpha)~in~B.$\\\\

\noindent{\it{\bf Theorem 2-20.}}  Let $B$ be a Banach $A-bimodule$ and $A$ has a $LBAI$. Suppose that $AB^{**},~B^{**}A\subseteq B$ and  every derivation $A^{**}$ into  $B^*$ is $weak^*-to-weak^*$ continuous. Then if $H^1(A,B^*)=0$, it follows that  $H^1(A^{**},B^{***})=0$.
\begin{proof} Take $B^{***}=B^*\oplus B^\perp$ where $B^\perp=\{b^{\prime\prime\prime}\in B^{***}:~b^{\prime\prime\prime}\mid_B=0\}$. Then we have
$$H^1(A^{**},B^{***})=H^1(A^{**},B^*)\oplus H^1(A^{**},B^\perp).$$
Since $H^1(A,B^*)=0$, it is similar to Theorem 2.17  that $H^1(A^{**},B^*)=0$. It suffices for the result to show that $H^1(A^{**},B^\perp)=0$.
 Without loss generality, let $(e_{\alpha})_{\alpha}\subseteq A$ be a $LBAI$ for $A$ such that  $e_{\alpha} \stackrel{w^*} {\rightarrow}e^{\prime\prime}$ in $A^{**}$ where $e^{\prime\prime}$ is a left unit for  $A^{**}$ with respect to the first Arens product. Let $a^{\prime\prime}\in A^{**}$ and suppose that $(a_{\beta})_{\beta}\subseteq A$ such that
 $a_{\beta} \stackrel{w^*} {\rightarrow}a^{\prime\prime}$ in $A^{**}$. Then, if we take $D\in Z^1(A^{**},B^\perp)$, for every $b^{\prime\prime}\in B^{**}$, by using the $weak^*-to-weak^*$ continuity  of $D$, we have
$$<D(a^{\prime\prime}), b^{\prime\prime}>= <D(e^{\prime\prime}a^{\prime\prime}), b^{\prime\prime}>=
\lim_\alpha\lim_\beta <(D(e_\alpha a_\beta),b^{\prime\prime}>$$$$
=\lim_\alpha\lim_\beta <(D(e_\alpha) a_\beta +e_\alpha D(a_\beta)),b^{\prime\prime}>
=\lim_\alpha\lim_\beta <D(e_\alpha) a_\beta ,b^{\prime\prime}>$$$$+\lim_\alpha\lim_\beta <e_\alpha D(a_\beta),b^{\prime\prime}>
=\lim_\alpha\lim_\beta <D(e_\alpha), a_\beta b^{\prime\prime}>$$$$+\lim_\alpha\lim_\beta < D(a_\beta),b^{\prime\prime}e_\alpha >=0.$$
It follows that $D=0$, and so the result is hold.\\ \end{proof}

\noindent{\it{\bf Corollary 2-21.}} Assume that $A$ is a Banach algebra with $LBAI$. Suppose that $A$ is two-sided ideal in $A^{**}$  and every derivation  $D:A^{**}\rightarrow A^{***}$ is $weak^*-to- weak^*$ continuous. Then if $A$ is weakly amenable, it follows that $A^{**}$ is weakly amenable.\\\\

Let $n\geq 0$ and suppose that $A$ is $(2n+1)-weakly$ amenable.  By condations of the proceeding corollary and by using [6, Corollary 1.14], we conclude that $A^{**}$ is weakly amenable.\\
Assume that $G$ is a locally compact group. We know that $L^1(G)$ is weakly amenable Banach algebra, see [15]. So by using proceeding corollary,  every $weak^*-to-weak^*$ continuous derivation from $L^1(G)^{**}$ into $L^1(G)^{***}$ is inner.\\

Let $A$ be a Banach Algebra. A dual Banach $A-bimodule$ $B$ is called normal if, for each $b\in B$ the map $a\rightarrow ax$ and $a\rightarrow xa$ from $A$ into $B$ is $weak^*-to-weak^*$ continuous.\\
A dual Banach algebra $A$ is Connes-amenable if, for every normal, dual Banach $A-bimodule$ $B$, every $weak^*-to-weak^*$ continuous derivation $D\in Z^1(A,B)$ is inner. Then we write $H^1_{w^*}(A,B)=0$.\\
A Banach algebra $A$ is called super-amenable if $H^1(A,B)=0$ for every Banach $A-bimodule$ $B$. It is clear that if $A$ is super-amenable, then $A$ is  amenable.\\
If $B$ is a Banach algebra and it is consisting a closed subalgebra $A$ such that $A$ is weakly dense in $B$. Then, it is easy to show that $A$ is a super-amenable if and only if $B$ is a super-amenable.\\
In the following, we will study some problems on the Connes-amenability of Banach algebras. In the  Theorem 2-19, with some conditions, we showed that $H^1(A^{**},B^{**})=0$. Now,  in the following theorem, by using some cocitions and super-amenability of Banach algebra $A$, we conclude that $H^1_{w^*}(A^{**},B^{**})=0$.\\
By Theorem 4.4.8 from [22], we know that if $A$ is Arens regular Banach algebra and $A$ is an ideal in $A^{**}$, then $A$ is amenable if and only if $A^{**}$ is Connes-amenable. Now in the following by using of super-amenability of Banach algebra $A$, we show that $A^{**}$ is Connes-amenable.\\
In the Theorem 2-19, for a Banach $A-bimodule$ $B$, with some conditions, we showed that $H^1(A^{**},B^{**})=0$. In the following, by using some new conditions, we show that $H_{w^*}^1(A^{**},B^{**})=0$ and as results in group algebra, for a  amenable compact group $G$, we show that $H^1_{w^*}(L^1(G)^{**},M(G)^{***})=0.$\\\\

\noindent{\it{\bf Theorem 2-22.}} Suppose that $A$ is supper-amenable  Banach algebra. Then we have the following statements.
\begin{enumerate}
\item $A^{**}$ is Connes-amenable Banach algebra with respect to the first Arens product.
\item  If $A$ has a $LBAI$ and for every Banach $A-bimodule$ $B$ with predual $X$, we have $AB^*,~B^*A\subseteq X$, then $H^1_{w^*}(A^{**},B^{**})=0$.
\end{enumerate}
\begin{proof}
\begin{enumerate}
\item Let  $B$ be a normal dual Banach $A-bimodule$ $B$ and let  be $D\in Z^1(A^{**},B^{**})$ be $weak^*-to-weak^*$ continuous derivation. Suppose that $\widetilde{D}=D\mid_A$. Then we have  $\widetilde{D}\in Z^1(A,B)$. Since $A$ is super-amenable, there is element  $b\in B$ such that $\widetilde{D}=\delta_b$. Suppose that $a^{\prime\prime}\in A^{**}$ and $(a_{\alpha})_{\alpha}\subseteq A$  such that  $a_{\alpha} \stackrel{w^*} {\rightarrow}a^{\prime\prime}$ in $A^{**}$. Since $\widetilde{D}$ is $weak^*-to-weak^*$ continuous derivation, we have
$$D(a^{\prime\prime})=D(w^*-\lim_\alpha a_\alpha)=w^*-\lim_\alpha D(a_\alpha)=w^*-\lim_\alpha \widetilde{D}(a_\alpha)$$$$=w^*-\lim_\alpha \delta_b(a_\alpha)=w^*-\lim_\alpha (a_\alpha b-ba_\alpha )=a^{\prime\prime}b-ba^{\prime\prime}.$$
It follows that $H^1_{w^*}(A^{**},B)=0$.

\item By using part (1), proof is similar to Theorem 2-19 (2).\\
\end{enumerate}
\end{proof}

\noindent{\it{\bf Theorem 2-23.}} Suppose that $A$ is  an amenable  Banach algebra with $LBAI$. If for every Banach $A-bimodule$ $B$, we have $AB^{**},~B^{**}A\subseteq B$, then $H^1_{w^*}(A^{**},B^{***})=0$.
\begin{proof} Proof is similar to the proceeding theorem.\\\end{proof}

Suppose that $G$ is a compact group. Then $L^1(G)$ is an ideal in $M(G)^{**}$. If $G$ is amenable, then by using proceeding theorem, we conclude the following equality.
$$H^1_{w^*}(L^1(G)^{**},M(G)^{***})=0.$$\\

\noindent{\it{\bf Corollary 2-24.}} Assume that $A$ is a weakly amenable Banach algebra with $LBAI$. Then if $A$ is an ideal in $A^{**}$, it follows that
$$H^1_{w^*}(A^{**},A^{***})=0.$$\\

\noindent{\it{\bf Example 2-25.}} Assume that $G$ is a compact group. Then we know that $L^1(G)$ has a bounded approximate identity [=$BAI$] and $L^1(G)$ is two-sided ideal in $L^1(G)^{**}$.  We know that $L^1(G)$ is weakly amenable. Hence it is clear $$H^1_{w^*}(L^\infty(G)^*,L^\infty(G)^{**})=0.$$\\

Let $B$ be a Banach $A-bimodule$. Then for every $b^{\prime\prime}\in B^{**}$, we define
$$L_{b^{\prime\prime}}(a^{\prime\prime})=b^{\prime\prime}a^{\prime\prime}~~~and~~~~
R_{b^{\prime\prime}}(a^{\prime\prime})=a^{\prime\prime}b^{\prime\prime},$$
for every $a^{\prime\prime}\in A^{**}$. These are the operation of left and right multiplication by $b^{\prime\prime}$ on $A^{**}$.\\
Let $B$ be a Banach $A-bimodule$. We say that $B$ is a left [resp. right] factors with respect to $A$, if $BA=B$ [resp. $AB=B$].\\

In the following, for super-amenable Banach algebra $A$ and a Banach $A-bimodule$ $C$, we have represented for every derivation $D$ from $A$ into $C$, and as a result in group algebras, for amenable locally compact group $G$, we show that every derivation $D$ from $L^1(G)$ into $L^\infty (G)$ is in the form $D=L_f$ where $f\in L^\infty (G)$.\\

\noindent{\it{\bf Theorem 2-26.}} Assume that $A$ is a super-amenable Banach algebra and $B$ is a Banach $A-bimodule$ such that $B$ factors on the left [resp. right]. Then there is a Banach $A-bimodule$ $(C,.)$ such that $C=B$ and
$Z^1(A,C)=\{L_{D^{\prime\prime}(e^{\prime\prime})}:~D\in Z^1(A,C)\}~~[resp.~~Z^1(A,C)=\{R_{D^{\prime\prime}(e^{\prime\prime})}:~D\in Z^1(A,C)\}]$ where $e^{\prime\prime}$ is a right [resp. left] identity for $A^{**}$.\\
\begin{proof} Certainly, every super-amenable Banach algebra is amenable. So, by using [22, Propostion 2.2.1], $A$ has a $BAI$ such as $(e_{\alpha})_{\alpha}$. Since $BA=B$, for every $b\in B$, there are $y\in B$ and $a\in A$ such that $b=ya$. Then we have
$$\lim_\alpha b e_{\alpha}=\lim_\alpha (ya) e_{\alpha}=\lim_\alpha y(a e_{\alpha})=ya=b.$$
It follows that $B$ has $RBAI$ as $(e_{\alpha})_{\alpha}\subseteq A$. Without loss generality, let $e^{\prime\prime}$ be a right unit for $A^{**}$ such that $e_{\alpha} \stackrel{w^*} {\rightarrow}e^{\prime\prime}$ in $A^{**}$.\\
Take $C=B$ and for every $a\in A$ and $x\in C$, we define $a.x=0$ and $x.a=xa$. It is clear that $(C,.)$ is a Banach $A-bimodule$. Suppose that $D\in Z^1(A,C)$. Then there is a $c\in C$ such that $D=\delta_c$. Then for every $a\in A$, we have
$$D(a)=\delta_c (a)=a.c-c.a=-ca.$$
Suppose that $x\in C$ and $x^\prime\in C^*$. Since $C.A=C$, there is $t\in C$ and $s\in A$ such that $x=t.s=ts$. Then we have
$$<x^\prime, xe_\alpha>=<x^\prime, tse_\alpha>=<x^\prime t,se_\alpha>\rightarrow
<x^\prime t,s>=<x^\prime, x>.$$
It follows that $xe_{\alpha} \stackrel{w} {\rightarrow}x$ in $C$. Since $D^{\prime\prime}$ is a $weak^*-to-weak^*$ continuous derivation, we have
$$D^{\prime\prime}(e^{\prime\prime})=D^{\prime\prime}(w^*-\lim_\alpha e_{\alpha})=w^*-\lim_\alpha D^{\prime\prime}(e_{\alpha})=w-\lim_\alpha D(e_{\alpha})$$$$=w-\lim_\alpha (ce_{\alpha})=-c.$$
Thus we conclude that $D(a)=D^{\prime\prime}(e^{\prime\prime})a$ for all $a\in A$. It follows that $$D=L_{D^{\prime\prime}(e^{\prime\prime})}.$$
On the other hand, since for every derivation $D\in Z^1(A,C)$, we have $L_{D^{\prime\prime}(e^{\prime\prime})}\in
Z^1(A,C)$, the result is hold.\\\end{proof}

\noindent{\it{\bf Corollary 2-27.}} Suppose that $A$ is  an amenable  Banach algebra and $B$ is a Banach $A-bimodule$ such that $AB^*=B^*$ [resp. $B^*A=B^*$]. Then there is a Banach $A-bimodule$ $(C,.)$ such that $C=B$ and
$Z^1(A,C^*)=\{L_{D^{\prime\prime}(e^{\prime\prime})}:~D\in Z^1(A,C^*)\}~~[resp.~~Z^1(A,C^*)=\{R_{D^{\prime\prime}(e^{\prime\prime})}:~D\in Z^1(A,C^*)\}]$ where $e^{\prime\prime}$ is a right [resp. left] identity for $A^{**}$.\\\\

\noindent{\it{\bf Corollary 2-28.}} Suppose that $A$ is  a super-amenable  Banach algebra and $B$ is a Banach $A-bimodule$. Then there is a Banach $A-bimodule$ $(C,.)$ such that $C=B$ and
$Z^1(A,C)=\{L_{c}:~c\in C\}~~[resp.~~Z^1(A,C)=\{R_{c}:~c\in C\}]$.\\\\

\noindent{\it{\bf Example 2-29.}}
\begin{enumerate}
\item  Let $G$ be an amenable locally compact group. Then there is a Banach $L^1(G)-bimodule$ such as $(L^\infty(G),.)$ such that $Z^1(L^1(G),L^\infty(G))=\{L_{f}:~f\in L^\infty(G)\}.$
\item  Let $G$ be locally compact group and $H^1(M(G),L^1(G))=0$. Then there is a Banach $M(G)-bimodule$ $(L^1(G),.)$ such that $Z^1(M(G),L^1(G))\\=\{L_{D^{\prime\prime}(e^{\prime\prime})}:~D\in Z^1(M(G),L^1(G))\}~~[or~~Z^1(M(G),L^1(G))=\{R_{D^{\prime\prime}(e^{\prime\prime})}:~\\D\in Z^1(M(G),L^1(G))\}]$ where $e^{\prime\prime}$ is a left [or right] identity for $L^1(G)^{**}$.\\\\

\end{enumerate}
\noindent{\it{\bf Theorem 2-30.}}  Let $B$ be a Banach $A-bimodule$ and suppose that $D:A\rightarrow B^*$ is a derivation. If $D^{\prime\prime}:A^{**}\rightarrow B^{***}$ is a derivation and $B^*\subseteq D^{\prime\prime}(A^{**})$, then $Z^\ell_{A^{**}}(B^{**})=B^{**}$.
\begin{proof} Since $D^{\prime\prime}:A^{**}\rightarrow B^{***}$ is a derivation, by [19, Theorem 4.2],
$D^{\prime\prime}(A^{**})B^{**}\subseteq B^*$. Due to $B^*\subseteq D^{\prime\prime}(A^{**})$, we have $B^*B^{**}\subseteq B^*$. Let   $(a_{\alpha}^{\prime\prime})_{\alpha}\subseteq A^{**}$ such that  $a^{\prime\prime}_{\alpha} \stackrel{w^*} {\rightarrow}a^{\prime\prime}$ in $A^{**}$. Assume that $b^{\prime\prime}\in B^{**}$. Then for every $b^\prime\in B^*$, since $b^\prime b^{\prime\prime}\in B^*$, we have
$$<b^{\prime\prime}a_{\alpha}^{\prime\prime},b^\prime>=<a_{\alpha}^{\prime\prime}b^\prime, b^{\prime\prime}>\rightarrow
<a^{\prime\prime},b^\prime b^{\prime\prime}>=<b^{\prime\prime}a^{\prime\prime},b^\prime >.$$
Thus $b^{\prime\prime}a^{\prime\prime}_{\alpha} \stackrel{w^*} {\rightarrow}b^{\prime\prime}a^{\prime\prime}$ in $B^{**}$, and so $b^{\prime\prime}\in Z^\ell_{A^{**}}(B^{**})$.\\ \end{proof}

\noindent{\it{\bf Corollary 2-31.}} Assume that $A$ is a Banach algebra and $D:A\rightarrow A^*$ is a derivation such that $A^*\subseteq D^{\prime\prime}(A^{**})$. Then if the linear mapping $D^{\prime\prime}:A^{**}\rightarrow A^{***}$
is a derivation, it follows that  $A$ is Arens regular.\\\\

\noindent{\it{\bf Lemma 2-32.}} Let $B$ be a Banach left $A-module$ and $B^{**}$ has a $LBAI$ with respect to $A^{**}$. Then $B^{**}$ has a left unit with respect to $A^{**}$.
\begin{proof}
Assume that $(e^{\prime\prime}_{\alpha})_{\alpha}\subseteq A^{**}$ is a $LBAI$ for $B^{**}$. By passing to a subnet, we may suppose that there is  $e^{\prime\prime}\in A^{**}$ such that $e^{\prime\prime}_{\alpha} \stackrel{w^*} {\rightarrow}e^{\prime\prime}$ in $A^{**}$. Then for every $b^{\prime\prime}\in B^{**}$ and $b^\prime\in B^*$, we have

$$<\pi_\ell^{***}(e^{\prime\prime},b^{\prime\prime}),b^\prime>=
<e^{\prime\prime},\pi_\ell^{**}(b^{\prime\prime},b^\prime)>=
\lim_\alpha <e_\alpha^{\prime\prime},\pi_\ell^{**}(b^{\prime\prime},b^\prime)>$$$$=
\lim_\alpha <\pi_\ell^{***}(e_\alpha^{\prime\prime},b^{\prime\prime}),b^\prime>=
<b^{\prime\prime},b^\prime>.$$
It follows that $\pi_\ell^{***}(e^{\prime\prime},b^{\prime\prime})=b^{\prime\prime}$.\\\end{proof}

\noindent{\it{\bf Corollary 2-33.}} Let $A$ be a Banach algebra and $A^{**}$ has a $LBAI$. Then $A^{**}$ has a left unit with respect to the first Arens product.\\\\

\noindent{\it{\bf Theorem 2-34.}} Let $A$ be a left strongly Arens irregular and suppose that $A^{**}$ is an amenable Banach algebra. Then we have the following assertions.
\begin{enumerate}
\item $A$ has an identity.

\item If $A$ is a dual Banach algebra, then $A$ is reflexive.
\end{enumerate}
\begin{proof}\begin{enumerate}
\item Since $A^{**}$ is amenable, it has a $BAI$.  By using the proceeding corollary, $A^{**}$ has an identity $e^{\prime\prime}$. So, the mapping $x^{\prime\prime}\rightarrow  e^{\prime\prime} x^{\prime\prime}=x^{\prime\prime}$ is $weak^*-to-weak^*$ continuous from $A^{**}$ into $A^{**}$. It follows that $e^{\prime\prime}\in Z_1(A^{**})=A$. Consequently, $A$ has an identity.

\item Assume that $E$ is predual Banach algebra for $A$. Then we have $A^{**}=A\oplus E^\bot$. Since $A^{**}$ is amenable, by using Theorem   , $A$ is amenable, and so $E^\bot$ is amenable. Thus $E^\bot$ has a $BAI$ such as $(e^{\prime\prime}_\alpha)_\alpha\subseteq E^\bot$. Since $E^\bot$ is a closed and $weak^*-closed$ subspace of $A^{**}$, without loss generality,  there is $e^{\prime\prime}\in E^\bot$ such that
   $$e^{\prime\prime}_\alpha\stackrel{w^*} {\rightarrow}e^{\prime\prime}~~and~~e^{\prime\prime}_\alpha\stackrel{\parallel~~\parallel} {\rightarrow}e^{\prime\prime}$$
 Then   $e^{\prime\prime}$ is a left identity for  $E^\bot$. On the other hand, for every $x^{\prime\prime}\in E^\bot$, since $E^\bot$ is an ideal in $A^{**}$, we have $x^{\prime\prime}e^{\prime\prime}\in E^\bot$. Thus for every $a^\prime\in A^*$, we have
 $$<x^{\prime\prime}e^{\prime\prime},a^\prime>=\lim_\alpha
 <(x^{\prime\prime}e^{\prime\prime})e^{\prime\prime}_\alpha,a^\prime>=\lim_\alpha
 <x^{\prime\prime}(e^{\prime\prime}e^{\prime\prime}_\alpha),a^\prime>
 $$$$=\lim_\alpha
 <x^{\prime\prime}e^{\prime\prime}_\alpha,a^\prime>= <x^{\prime\prime},a^\prime>.$$
 It follows that $x^{\prime\prime}e^{\prime\prime}=x^{\prime\prime}$, and so  $e^{\prime\prime}$ is a right identity for $E^\bot$. Consequently, $e^{\prime\prime}$ is a two-sided identity for $E^\bot$.
Now, let $a^{\prime\prime}\in A^{**}$. Then we have the following equalities.
$$e^{\prime\prime}a^{\prime\prime}= (e^{\prime\prime}a^{\prime\prime})e^{\prime\prime}=
 e^{\prime\prime}(a^{\prime\prime}e^{\prime\prime})=a^{\prime\prime}e^{\prime\prime}.$$
Thus we have $e^{\prime\prime}\in Z_1(A^{**})=A$. It follows that  $e^{\prime\prime}=0$, and so $E^\bot=0$. Consequently, we have $A^{**}=A$.\\\end{enumerate}\end{proof}

 Let $G$ be a locally compact group. Then if $L^1(G)^{**}$ or $M(G)^{**}$ are amenable, by using proceeding theorem part (1) and (2), respectively,  we conclude that $G$ is finite group.\\\\
\newpage

\section{ \bf  Weak amenability of Banach algebras  }

\noindent In this section,  for Banach $A-module$ $B$, we introduce some new concepts as  $left-weak^*-to-weak$ convergence property [ $=Lw^*wc-$property] and $right-weak^*-to-weak$ convergence property [ $=Rw^*wc-$property] with respect to $A$ and we show that if $A^*$ and $A^{**}$, respectively,  have $Rw^*wc-$property and $Lw^*wc-$property and  $A^{**}$ is weakly amenable, then $A$ is weakly amenable. We also show the relations between a derivation $D:A\rightarrow A^*$ and this new concepts.
Now in the following, for left and right Banach $A-module$ $B$, we define, respectively, $Lw^*wc-$property and $Rw^*wc-$property concepts with some examples. \\

\noindent{\it{\bf Definition 3-1.}} Assume that  $B$ is a left Banach $A-bimodule$. We say that $b^\prime \in B^*$ has $left-weak^*-to-weak$ convergence property [ $=Lw^*wc-$property] with respect to $A$, if for every bounded net $(a_\alpha)_\alpha \subseteq A$, $b^\prime a_\alpha  \stackrel{w^*} {\rightarrow}0$ implies
 $b^\prime a_\alpha  \stackrel{w} {\rightarrow}0$.\\
 When every $b^\prime \in B^*$ has   $Lw^*wc$-property with respect to $A$, we say that $B^*$ has $Lw^*wc-$property with respect to $A$.\\
 The definition of $right-weak^*-to-weak$ convergence property [$=Rw^*wc-$property] with respect to $A$ is similar and if $b^\prime \in B^*$ has $left-weak^*-to-weak$ convergence property and $right-weak^*-to-weak$ convergence property, then we say that $b^\prime \in B^*$ has $weak^*-to-weak$ convergence property  [$=w^*wc-$property].\\\\\\
\noindent{\it{\bf Example 3-2 .}} \begin{enumerate}

 \item ~ Every reflexive Banach $A-module$ has $w^*wc-$property.

  \item ~Let $\Omega$ be a compact group and suppose that $A=C(\Omega)$ and $B=M(\Omega)$. Let $(a_\alpha)_\alpha \subseteq A$ and $\mu\in B$. Suppose that $\mu * a_\alpha  \stackrel{w^*} {\rightarrow}0$, then for each $a\in A$, we have
 $$<\mu * a_\alpha ,a>=<\mu , a_\alpha *a>=\int_\Omega (a_\alpha *a) d\mu\rightarrow 0.$$
 We set $a=1_\Omega$ . Then $\mu(a_\alpha)\rightarrow 0$. Now let $b^\prime\in B^*$. Then
 $$<b^\prime, \mu * a_\alpha>=<a_\alpha b^\prime, \mu  >=\int_\Omega a_\alpha b^\prime d\mu\leq \parallel b^\prime \parallel~\mid\int_\Omega a_\alpha d\mu \mid=\parallel b^\prime \parallel~\mid\mu( a_\alpha ) \mid\rightarrow 0.$$
It follows that  $\mu * a_\alpha  \stackrel{w} {\rightarrow}0$, and so that $\mu$ has $Rw^*wc-$property with respect to $A$.\\

 \end{enumerate}

 \noindent{\it{\bf Theorem 3-3.}} Let $A$ be a Banach algebra and suppose that $A^*$ and $A^{**}$, respectively,  have $Rw^*wc-$property and $Lw^*wc-$property with respect to $A$. If $A^{**}$ is weakly amenable, then $A$ is weakly amenable.
\begin{proof} Assume that $a^{\prime\prime}\in A^{**}$ and $(a_\alpha)_\alpha \subseteq A$ such that $a_\alpha \stackrel{w^*} {\rightarrow}a^{\prime\prime}$. Then for each $a^\prime\in A^*$, we have $a_\alpha a^\prime\stackrel{w^*} {\rightarrow}a^{\prime\prime}a^\prime$ in $A^*$. Since $A^*$ has $Rw^*wc-$property with respect to $A$, $a_\alpha a^\prime\stackrel{w} {\rightarrow}a^{\prime\prime}a^\prime$ in $A^*$. Then for every $x^{\prime\prime}\in A^{**}$, we have
$$<x^{\prime\prime}a_\alpha , a^\prime >=<x^{\prime\prime},a_\alpha  a^\prime >\rightarrow
<x^{\prime\prime},a^{\prime\prime}  a^\prime >=<x^{\prime\prime}a^{\prime\prime},  a^\prime >.$$
It follows that $x^{\prime\prime}a_\alpha \stackrel{w^*} {\rightarrow}x^{\prime\prime}a^{\prime\prime}$. Since $A^{**}$ has $Lw^*wc-$property with respect to $A$, $x^{\prime\prime}a_\alpha \stackrel{w} {\rightarrow}x^{\prime\prime}a^{\prime\prime}$. If $D:A\rightarrow A^*$ is a bounded derivation, we extend it to a bounded linear mapping  $D^{\prime\prime}$ from $A^{**}$ into $A^{***}$. Suppose that  $a^{\prime\prime}, b^{\prime\prime}\in A^{**}$ and $(a_\alpha)_\alpha , (b_\beta)_\beta\subseteq A$ such that $a_\alpha \stackrel{w^*} {\rightarrow}a^{\prime\prime}$ and $b_\beta \stackrel{w^*} {\rightarrow}b^{\prime\prime}$. Since  $x^{\prime\prime}a_\alpha \stackrel{w} {\rightarrow}x^{\prime\prime}a^{\prime\prime}$ for every $x^{\prime\prime}\in A^{**}$, we have
$$\lim_\alpha <D^{\prime\prime}(b^{\prime\prime}), x^{\prime\prime}a_\alpha>=<D^{\prime\prime}(b^{\prime\prime}), x^{\prime\prime}a^{\prime\prime}>.$$
In the following we take limit on the $weak^*$ topologies.  Thus we have
$$\lim_\alpha \lim_\beta D(a_\alpha)b_\beta=D^{\prime\prime}(a^{\prime\prime})b^{\prime\prime}.$$
Consequently, it follows that
$$D^{\prime\prime}(a^{\prime\prime}b^{\prime\prime})=\lim_\alpha \lim_\beta D(a_\alpha b_\beta)=\lim_\alpha \lim_\beta D(a_\alpha)b_\beta +\lim_\alpha \lim_\beta a_\alpha D(b_\beta)$$$$=D^{\prime\prime}(a^{\prime\prime})b^{\prime\prime}+ a^{\prime\prime}D^{\prime\prime}(b^{\prime\prime}).$$
Since $A^{**}$ is weakly amenable, there is $a^{\prime\prime\prime}\in A^{***}$ such that $D^{\prime\prime}=\delta_{a^{\prime\prime\prime}}$. Then we have
$D=D^{\prime\prime}\mid_{A^*}=\delta_{a^{\prime\prime\prime}}\mid_{A^{*}}$. Hence for each $x^\prime\in A ^*$, we have $D=x^\prime a^{\prime\prime\prime}-a^{\prime\prime\prime}x^\prime$. Take $a=a^{\prime\prime\prime}\mid_{A^* }$. It follows that $H^1(A,A^*)=0$.\\
\end{proof}

\noindent{\it{\bf Theorem 3-4.}} Let $A$ be a Banach algebra and suppose that $D:A\rightarrow A^*$ is a surjective  derivation. If $D^{\prime\prime}$ is a derivation, then we have the following assertions.
\begin{enumerate}
\item ~$A^*$ and $A^{**}$, respectively,  have $w^*wc-$property and $Lw^*wc-$property with respect to $A$.
\item ~For every $a^{\prime\prime}\in A^{**}$, the mapping $x^{\prime\prime}\rightarrow a^{\prime\prime}x^{\prime\prime}$ from $A^{**}$ into $A^{**}$ is $weak^*-to-weak$ continuous.
\item ~ $A$ is Arens regular.
\item ~ If $A$ has $LBAI$, then $A$ is reflexive.\\
 \end{enumerate}

\begin{proof}
\begin{enumerate}
\item ~Since $D$ is surjective, $D^{\prime\prime}$ is surjective, and so by using [19, Theorem 2.2], we have $A^{***}A^{**}\subseteq D^{\prime\prime}(A^{**})A^{**}\subseteq A^*$. Suppose that $a^{\prime\prime}\in A^{**}$ and  $(a_\alpha)_\alpha \subseteq A$ such that
    $a_\alpha \stackrel{w^*} {\rightarrow}a^{\prime\prime}$. Then for each   $x^\prime\in A^*$, we have  $x^\prime a_\alpha \stackrel{w^*} {\rightarrow}x^\prime a^{\prime\prime}$. Since $A^{***}A^{**}\subseteq A^*$, $x^\prime a^{\prime\prime}\in A^*$. Then for every  $x^{\prime\prime}\in A^{**}$, we have
$$<x^{\prime\prime}, x^\prime a_\alpha>=<x^{\prime\prime} x^\prime , a_\alpha>\rightarrow   <a^{\prime\prime} , x^{\prime\prime}x^\prime >=<x^\prime a^{\prime\prime} , x^{\prime\prime} >=< x^{\prime\prime},x^\prime a^{\prime\prime}  >.$$
 It follows that  $x^\prime a_\alpha \stackrel{w} {\rightarrow}x^\prime a^{\prime\prime}$ in $A^*$. Thus by easy calculation, we conclude that $x^\prime$ has  $Lw^*wc-$property with respect to $A$. It is to similar that $x^\prime$ has  $Rw^*wc-$property with respect to $A$, and so $A^*$ has $w^*wc-$property. \\
For next part, suppose that $x ^{\prime\prime\prime} \in A^{***}$.  Since $A^{***}A^{**}\subseteq A^*$,   $x^{\prime\prime}a_\alpha \stackrel{w^*} {\rightarrow}x^{\prime\prime}a^{\prime\prime}$ for each $x^{\prime\prime}\in A^{**}$. Then
$$<x ^{\prime\prime\prime}, x^{\prime\prime}a_\alpha>=<x ^{\prime\prime\prime} x^{\prime\prime}, a_\alpha>\rightarrow <x ^{\prime\prime\prime} x^{\prime\prime}, a^{\prime\prime}>=<x ^{\prime\prime\prime}, x^{\prime\prime} a^{\prime\prime}>.$$
It follows that $x^{\prime\prime}a_\alpha \stackrel{w} {\rightarrow}x^{\prime\prime}a^{\prime\prime}$.  Thus by easy calculation, we conclude that $x^{\prime\prime}$ has  $Lw^*wc-$property with respect to $A$.
\item Suppose that $(a^{\prime\prime}_\alpha)_\alpha \subseteq A^{**}$ and  $a^{\prime\prime}_\alpha \stackrel{w^*} {\rightarrow}a^{\prime\prime}$. Let $x^{\prime\prime}\in A^{**}$. Then for every $x ^{\prime\prime\prime} \in A^{***}$, since $A^{***}A^{**}\subseteq A^*$, we have
$$<x ^{\prime\prime\prime}, x^{\prime\prime}a^{\prime\prime}_\alpha>=<x ^{\prime\prime\prime} x^{\prime\prime}, a^{\prime\prime}_\alpha>\rightarrow <x ^{\prime\prime\prime} x^{\prime\prime}, a^{\prime\prime}>=
<x ^{\prime\prime\prime}, x^{\prime\prime} a^{\prime\prime}>.$$

\item ~~It  follows from (2).
\item ~Let $(e_\alpha)_\alpha \subseteq A$ be a $BLAI$ for $A$. Then without loss generality, let $e^{\prime\prime}$ be a left unit for $A^{**}$ such that  $e_\alpha \stackrel{w^*} {\rightarrow}e^{\prime\prime}$. Suppose that $(a^{\prime\prime}_\alpha)_\alpha \subseteq A^{**}$ and  $a^{\prime\prime}_\alpha \stackrel{w^*} {\rightarrow}a^{\prime\prime}$. Then for every $a ^{\prime\prime\prime} \in A^{***}$, since $A^{***}A^{**}\subseteq A^*$, we have
$$<a ^{\prime\prime\prime},a^{\prime\prime}_\alpha>= <a ^{\prime\prime\prime},e^{\prime\prime}a^{\prime\prime}_\alpha>=<a ^{\prime\prime\prime}e^{\prime\prime}, a^{\prime\prime}_\alpha>\rightarrow <a ^{\prime\prime\prime}e^{\prime\prime}, a^{\prime\prime}>=<a ^{\prime\prime\prime}, a^{\prime\prime}>.$$
 It follows that   $a^{\prime\prime}_\alpha \stackrel{w} {\rightarrow}a^{\prime\prime}$. Consequently $A$ is reflexive.\\

 \end{enumerate}
\end{proof}

\noindent{\it{\bf Corollary 3-5.}} Let $A$ be a Banach algebra and suppose that $D:A\rightarrow A^*$ is a surjective  derivation. The following statement are equivalent.
\begin{enumerate}
\item ~$A^*$ and $A^{**}$, respectively,  have $Rw^*wc-$property and $Lw^*wc-$property with respect to $A$.
\item ~For every $a^{\prime\prime}\in A^{**}$, the mapping $x^{\prime\prime}\rightarrow a^{\prime\prime}x^{\prime\prime}$ from $A^{**}$ into $A^{**}$ is $weak^*-to-weak$ continuous.\\\\
 \end{enumerate}

\noindent{\it{\bf Problems.}} \\
i) By notice to Theorem 2-19, dose the following assertions hold.
\begin{enumerate}
\item For compact group $G$, we have $H^1(L^1(G)^{**},M(G))=0?$
\item If $c_0$ is weakly amenable, $H^1(\ell^\infty,\ell^1)=0?$
\end{enumerate}
ii) Suppose that  $S$ is a compact semigroup. Dose $L^1(S)^*$ and $M(S)^*$ have   $Lw^*wc-$property or $Rw^*wc-$property?\\\\

\bibliographystyle{amsplain}

\end{document}